\documentclass{article}

\usepackage{graphicx}

\title{Seshadri constants on ruled surfaces: the rational and the elliptic cases.}
\author{Luis Fuentes Garc\'{\i}a}
\date{}

\newtheorem{teo}{Theorem}[section]
\newtheorem{defin}[teo]{Definition}
\newtheorem{prop}[teo]{Proposition}
\newtheorem{cor}[teo]{Corollary}
\newtheorem{lemma}[teo]{Lemma}
\newtheorem{example}[teo]{Example}
\newtheorem{rem}[teo]{Remark}

\font\euf=eufm10 at 12pt

\def\dfrac#1#2{\displaystyle{{#1}\over{#2}}}

\def\g2{\pi}

\def\e{\mbox{\euf e}}

\def\P{{\bf P}}

\newcommand\E{{\cal E}}
\newcommand\Te{{\cal O}}
\newcommand\C{{C}}

\def\ZZ{\leavevmode\hbox{$\rm Z$}}
\def\QQ{\leavevmode\hbox{$\rm I\!\!\!Q$}}
\def\RR{\leavevmode\hbox{$\rm I\!R$}}

\def\impp{{\quad\Rightarrow\quad}}

\def\qed{\hspace{\fill}$\rule{2mm}{2mm}$}

\newcommand\lrw{\longrightarrow}

\begin{document}

\maketitle

\vspace{0.1cm}

\begin{abstract}

We study the Seshadri constants on geometrically ruled surfaces. The unstable case is completely solved. Moreover, we give some bounds for the stable case. We apply these results to compute the Seshadri constant of the rational and elliptic ruled surfaces. Both cases are completely determined. The elliptic case provides an interesting picture of how particular is the behavior of the Seshadri constants.



\end{abstract}

\section{Introduction.}

The Seshadri constants were introduced by Demailly in \cite{De92}. They measure the local positivity of 
ample line bundles on algebraic varieties. They can be useful to obtain results about the global generation, very ampleness or $s$-jets generation of adjoint linear systems (see \cite{EiKuLa95}).

However, in general, they are very hard to compute explicitly. Even for surfaces, there are not many examples where the Seshadri constants are known. They were computed for simple abelian surfaces by Th. Bauer and T. Szemberg (see \cite{Ba98}). Recently, Ch. Schultz gave explicit values for Seshadri constants on products of two elliptic curves (see \cite{Sc04}). General bounds for Seshadri constants on surfaces are given in \cite{Ba99}, \cite{Na04} or \cite{St98}.

An interesting open question is if they are always rational. There are not known examples where the Seshadri constants take irrational values. Note that, for surfaces with Picard number great than $1$, the rationality of the Seshadri constant implies the existence of curves with high multiplicity at a point. These curves will be called Seshadri exceptional curves (see \cite{EiKuLa95}).

Here, we study the Seshadri constants on geometrically ruled surfaces. We give explicit values when the ruled surface corresponds to an unstable vector bundle (see Theorem \ref{inestables} and Corollary \ref{trivial}). These results provide immediately the Seshadri constants on the rational ruled surfaces.

However, the most interesting cases appear when the ruled surface corresponds to a stable (or semi-stable) vector bundle.  A useful nice fact will be that the Seshadri exceptional curves are preserved by elementary transformations. We will give bounds for the Seshadri constants on the stable (or semi-stable) ruled surfaces (see Theorem \ref{principal}). Moreover, we will give explicit values of the Seshadri constants for particular semi-stable ruled surfaces. As an example of the interest  of the Seshadri constants on these cases, we will prove the following theorem:

\begin{teo}

Given a real number $\delta>0$ and a smooth curve $X$ of genus $>0$, there is a stable ruled surface $S$, an ample divisor $A$ on $S$ and a point $x\in S$ such that
$$
\sqrt{A^2}-\delta<\epsilon(A,x)<\sqrt{A^2}.
$$

\end{teo}

The last section of the article is dedicated to compute the Seshadri constants on the elliptic ruled surfaces. We will give explicitly the Seshadri constants of any ample divisor at any point of an arbitrary elliptic ruled surface. In particular, we show that they are always rational. The unstable cases will be computed with the general theory. The stable cases will be solved by using involutions of the elliptic ruled surfaces. They will allow us construct divisors with high multiplicity. 

The results provide a very nice picture of how particular is the behavior of the Seshadri constants. For example, if $S_{-1}$ is the indecomposable elliptic ruled surface with invariant $e=-1$, we obtain:

\begin{teo}

Let $A\equiv aX_0+bf$ be an ample divisor on $S_{-1}$. Let $x\in S_{-1}$.

\begin{enumerate}

\item If $x\in T_k$ with $k>1$ and 

\begin{enumerate}

\item $\frac{a}{2(b+\frac{1}{2}a)}\leq 1$ then
$$
\epsilon(A,x)=a.
$$

\item $1<\frac{a}{2(b+\frac{1}{2}a)}\leq \frac{k^2}{(k-2)^2}$ when $k$ is even or $1<\frac{a}{2(b+\frac{1}{2}a)}\leq \frac{k^2+1}{k^2-4k+5}$ when $k$ is odd, then

$$
\epsilon(A,x)=\frac{(k-1)a+kb}{k-1}.
$$

\item $\frac{a}{2(b+\frac{1}{2}a)}\geq \frac{k^2}{(k-2)^2}$ when $k$ is even or $\frac{a}{2(b+\frac{1}{2}a)}\geq \frac{k^2+1}{k^2-4k+5}$ when $k$ is odd, then

$$
\epsilon(A,x)=\frac{2n(na+(n+1)b)}{2n^2-1}
\hbox{ with }
n=\left[\frac{1}{\sqrt{\frac{a}{2(b+\frac{1}{2}a)}}-1}\right]+1.
$$

\end{enumerate}

\item If $x\not\in T$ and 
\begin{enumerate}

\item $\frac{a}{2(b+\frac{1}{2}a)}\leq 1$ then
$$
\epsilon(A,x)=a.
$$

\item $\frac{a}{2(b+\frac{1}{2}a)}> 1$ then
$$
\epsilon(A,x)=\frac{2n(na+(n+1)b)}{2n^2-1}
\hbox{ with }
n=\left[\frac{1}{\sqrt{\frac{a}{2(b+\frac{1}{2}a)}}-1}\right]+1.
$$

\end{enumerate}

\end{enumerate}

\end{teo}

Note that, for any $k>1$, the set $T_k$ is an effective divisor and $T=\bigcup_{k>1}T_k$ is a countable union of such divisors (see Section \ref{elipticas}). If $x\in S_{-1}/T$ ($x$ is a very general point) then the Seshadri constants do not depend on $x$. However, at points in the curves $T_k$ they have a special behavior. This justifies that most of the general results about Seshadri constants on surfaces are given for very general points.

We refer to \cite{FuPe} and \cite{Ha} for a systematic
development of the projective theory of ruled surfaces and to
\cite{Ba99} for the main properties of the Seshadri constants on surfaces. 

I  would like to express my hearty thanks to Professor Manuel Pedreira for valuable discussions and his warm encouragement.

\section{Preliminaries on Seshadri constants.}

Let $S$ be a smooth surface. We will denote by $Nef(S)$ the cone of nef $\RR$-divisors on $S$ and by $NE(S)$ the convex cone generated by reduced and irreducible curves on $S$.

 Let $f:\tilde{S}\lrw S$ be the the blowing up of $S$ at $x\in S$ and $E$ the exceptional divisor. 
Let $A$ be a nef divisor on $S$. We define de Seshadri constant of $A$ at $x$ as:
$$
\epsilon(A,x)=sup\{\epsilon\in \RR |\,f^*A-\epsilon E\hbox{ is nef }\}.
$$
Equivalently:
$$
\epsilon(A,x)=inf_{C\ni x}\left\{{\frac{A\cdot C}{mult_x(C)}|\,\hbox{ $C$ irreducible curve through $x$}}\right\}.
$$
A well known general bound for the Seshadri constant is:
$$
\epsilon(A,x)\leq \sqrt{A^2}.
$$
Furthermore, when the constant does not reach this upper bound there is a curve $C$ such that $\epsilon(A,x)=\frac{A\cdot C}{mult_x (C)}$ (see \cite{Ba99} and \cite{St98} ). This curve verifies:
$$
C^2<mult_x(C)^2.
$$

\begin{defin}

An irreducible curve $C$ passing through $x$ with multiplicity $m\geq 1$ and $C^2<m^2$ is called a Seshadri exceptional curve based at $x$.

\end{defin}

Given a Seshadri exceptional curve $C$based at $x$, we can define the continuous function (see \cite{Sc04}):
$$
q_C:Nef(X)\lrw \RR, \quad q_C(A,x)=\frac{A\cdot C}{mult_x (C)}.
$$
With this notation:
$$
\epsilon(A,x)=min(\{q_C(A,x)|\,\hbox{$C$ is Seshadri exceptional based at $x$}\}\cup \{\sqrt{A^2}\}).
$$

\begin{defin}

Let $C$ be a Seshadri exceptional curve based at $x$. The open set of nef divisors $A$ verifying:
$$
q_C(A,x)<\sqrt{A^2}
$$
is called the influence area of $C$  and it will be denoted by $Q_C$.

\end{defin}

If $C$ is a nef Seshadri exceptional curve based at $x$, it is clear that $C\in Q_C$. In fact, $C$ is the unique Seshadri exceptional curve based at $x$ in $\overline{Q_C}$:

\begin{lemma}\label{restriccion}

If $C,D$ are different Seshadri exceptional curves based at $x$, then $D\not\in \overline{Q_C}$.

\end{lemma}
{\bf Proof:} If $D$ is Seshadri exceptional and $D\in \overline{Q_C}$ then, $$\frac{D\cdot C}{mult_x(C)}\leq \sqrt{D^2}\impp D\cdot C\leq mult_x(C)\sqrt{D^2}<mult_x(C)\,mult_x(D).$$ Since $D$ are $C$ are different curves this is impossible. \qed

\begin{cor}

If $C$ is a nef Seshadri exceptional curve based at at $x$, then $\epsilon(C,x)=q_C(C,x)=\frac{C^2}{mult_x(C)}$. \qed

\end{cor}

It is clear that if we want to compute the Seshadri constant of any nef divisor at a point $x$ we must know all the Seshadri exceptional curves based at $x$ and their influence areas. In particular, with this information we can know when irrational Seshadri constants can appear:

\begin{lemma}

Let $S$ be a smooth surface with Picard number $\rho(S)>1$ and let $x\in S$. The Seshadri constant at $x$ is rational for all nef divisors if and only if 
$$
Nef(S)=\overline{\bigcup\{Q_C|\,\hbox{$C$ is Seshadri exceptional based at $x$}\}}.
$$

\end{lemma}
{\bf Proof:} If the equality holds, then for any nef divisor $A$ there is a Seshadri exceptional curve $C$ based at $x$ such that $q_C(A,x)\leq \sqrt{A^2}$, so $\epsilon(A,x)$ is rational. On the other hand, if the equaility does not hold, then there is an open $U$ set in $Nef(X)$, such that $\epsilon(A,x)=\sqrt{A^2}$ for any $A\in U$. If $\rho(S)>1$ we can construct a divisor $A\in U$ with $\sqrt{A^2}$ irrational. \qed

\section{Preliminaries on ruled surfaces.}

A {\it geometrically ruled surface}, or simply a {\it
ruled surface}, will be a $\P^1$-bundle over a smooth curve $X$ of genus $g>0$. It will be
denoted by $\pi: S=\P(\E_0)\lrw X$ with fibre $f$. We will suppose that
$\E_0$ is a normalized sheaf and $X_0$ is the section of minimum self-intersection
that corresponds to the surjection $\E_0\lrw \Te_X(\e)\lrw 0$, with $\bigwedge^2\E\cong
\Te_X(\e)$ and  $e=-deg(\e)$. We know that  $Num(S)=\ZZ X_0\oplus \ZZ f$ (see \cite{Ha},V,2 and \cite{FuPe}).

When $\E_0\equiv \Te_X\oplus \Te_X(\e)$ we say that the ruled surface is decomposable. In this case we will denote by $X_1$ the irreducible curves of the linear system $|X_0-\e f|$. They does not meet $X_0$ and have self-intersection $e$.

A curve $C\subset S$ will be said $n$-secant when $C\cdot f=n$. In particular, if $C\cdot f=1$ we say that $C$ is unisecant.

We will use the following results which characterize the nef and ample divisors on ruled surfaces:

\begin{prop}\label{noestable}

Let $S$ be a ruled surface with invariant $e>0$:

\begin{enumerate}

\item If $Y\equiv aX_0+bf$ is an irreducible curve $\neq X_0,f$ then $a>0$ and $b>ae$.

\item A divisor $A\equiv aX_0+bf$ is nef if and only if $a\geq 0$ and $b\geq ae$.

\item A divisor $A\equiv aX_0+bf$ is ample is and only if $a>0$ and $b> ae$.

\end{enumerate}

\end{prop}

\begin{prop}\label{estable}

Let $S$ be a ruled surface with invariant $e\leq 0$:

\begin{enumerate}

\item If $C\equiv aX_0+bf$ is an irreducible curve  $\neq X_0,f$ then  either $a=1$ and $b\geq 0$ or $a>1$ and $b\geq \frac{1}{2}ae$.

\item A divisor $A\equiv aX_0+bf$ is nef if and only if $a\geq 0$ and $b\geq \frac{1}{2}ae$.

\item A divisor $A\equiv aX_0+bf$ is ample if and only if $a>0$ and $b> \frac{1}{2}ae$.

\end{enumerate}

\end{prop}

\begin{rem}

If $S$ is a ruled surface with invariant $e\leq 0$, the convex cone of effective curves is generated by the $\QQ$-divisors $Z\equiv X_0+\frac{1}{2}ef$ and $f$. With this notation, a divisor $C\equiv aZ+b'f$ is nef if and only if $a\geq 0$ and $b'\geq 0$. Moreover the intersection pairing is determined by $f^2=0$, $Z^2=0$ and $f\cdot Z=1$. The advantage of these generators is that they does not depend on $e$. Note that $aX_0+bf\equiv aZ+(b-\frac{1}{2}ae)f$. In order to respect the classical notation we will consider the generators $f$ and $X_0$, but we will express the results in terms of $a$ and $(b-\frac{1}{2}ae)$.

\end{rem}

\section{Seshadri constants on ruled surfaces.}

\subsection{Seshadri exceptional curves on ruled surfaces.}

In this section we study the existence of Seshadri exceptional curves on a ruled surface $S$ with invariant $e$. It is clear that the fibres $f$ are Seshadri exceptional curves at any point. Moreover, we can compute their influence area. If $A\equiv aX_0+bf$ is a nef divisor then,
$$
A\cdot f<\sqrt{A^2}\iff a<\sqrt{2a(b-\frac{1}{2}ae)} \iff b-\frac{1}{2}ae> \frac{1}{2}a\quad\hbox{ and }\quad a>0.
$$
By Lemma \ref{restriccion}, this gives a first restriction for the rank where we can hope to find Seshadri exceptional curves. From this, we obtain some results:

\begin{lemma}\label{fibras}

The fibres are Seshadri exceptional curves based at any point $x\in S$. Moreover, their influence area is:
$$
Q_f=\{aX_0+bf\in Nef(S)|\, \frac{a}{2(b-\frac{1}{2}ae)}<1,\, a>0\}.
$$

\end{lemma} \qed

\begin{cor}\label{rank1}

Let $C\equiv aX_0+bf$ with $a>0$ be a Seshadri exceptional curve on $S$. Then $$\frac{a}{2(b-\frac{1}{2}ae)}> 1.$$
\end{cor} \qed

\begin{lemma}\label{unisecantes}

Let $C\equiv X_0+bf$ be an unisecant Seshadri exceptional curve based at $x$. Then $C\equiv X_0$ and $e\geq 0$. In particular, if $e=0$ and $C\equiv X_0$:
$$
Q_{C}=\{aX_0+bf\in Nef(S)|\, \frac{a}{2b}>\frac{1}{4}\}.
$$
\end{lemma}
{\bf Proof:} By the previous Corollary we know that $b-\frac{1}{2}e<\frac{1}{2}$. Applying the Propositions \ref{noestable} and \ref{estable} we see that the unique possibility is $e\geq 0$ and $b=0$. In this case, $C^2=0$ so $C$ is a Seshadri exceptional curve and we can compute its influence area. \qed

\begin{cor}\label{unicas}

Suppose that $e=0$ and let $C\equiv X_0$ be an irreducible curve. Let $x$ be a point of $C$. The unique Seshadri exceptional curves  based at $x$ are $C$ and $f$.

\end{cor}
{\bf Proof:} By the Lemma \ref{fibras} and the Lemma \ref{unisecantes}, $Q_f\cup Q_{C}=Nef(S)=\overline{NE}(S)$. Now the conclusion follows from Lemma \ref{restriccion}. \qed

Now, we give a bound for the multiplicity of an irreducible curve at any point. We will use that an $a$-secant curve on $S$, defines a $a:1$ map between $C$ and the base curve $X$.

\begin{lemma}\label{cota}

Let $C\equiv aX_0+bf$ an irreducible curve on $S$ with $a\geq 2$. Let $x\in C$ and $m=mult_x(C)$. Then:
$$
m(m-1)\leq C^2-2(b-\frac{1}{2}ae) 
$$
or equivalently
$$
m(m-1)\leq 2(a-1)(b-\frac{1}{2}ae).
$$

\end{lemma}
{\bf Proof:} Let $\tilde{C}$ be the normalization of $C$. Since $C$ meets each generator at $a$ points, we have and induced map $\gamma:\tilde{C}\lrw X$ of degree $a$. By the Hurwitz formula:
$$
2g(\tilde{C})-2=a(2g-2)+deg R\impp 2g(\tilde{C})-2-a(2g-2)\geq 0
$$
where $g$ and $g(\tilde{C})$ are respectively the genus of $X$ and $\tilde{C}$, and $R$ is the (effective) ramification divisor. Moreover, $g(\tilde{C})\leq p_a(C)-\frac{1}{2}m(m-1)$ where $p_a(C)$ is the arithmetic genus of $C$. From this:
$$
2p_a(C)-2-a(2g-2)\geq m(m-1).
$$
On the other hand, by the adjunction formula, we have $2p_a(C)-2=C\cdot (K_S+C)$. Combining this fact with the  inequality above:
$$
C\cdot (K_S+C)-a(2g-2)\geq m(m-1).
$$
The canonical divisor of the ruled surface is $K_S\equiv -2X_0+(2g-2-e)f$, so 
$$
C\cdot (K_S+C)-a(2g-2)=C^2+2ae+a(2g-2)-ae-2b-a(2g-2)=C^2-2(b-\frac{1}{2}ae).
$$ \qed

\begin{cor}

If $C$ is a Seshadri exceptional curve based at $x$ with $C^2>0$ then $\sqrt{C^2}$ is irrational and $mult_x(C)=[\sqrt{C^2}]+1$.

\end{cor}
{\bf Proof:} Let $m=mult_x(C)$. Since $C$ is Seshadri exceptional, $C^2<m^2$. On the other hand, since $C^2>0$, the curve $C$ is not a fibre, and by the Corollary \ref{unisecantes} it is not an unisecant curve. Thus we can apply the previous theorem. We obtain:
$$
m(m-1)<C^2<m^2 \impp (m-1)^2<C^2<m^2 \impp m-1<\sqrt{C^2}<m.
$$
Since $m$ is an integer number, $m=[\sqrt{C^2}]+1$ and $\sqrt{\C^2}$ is irrational. \qed

\begin{rem}

The same result is also true for an arbitrary surface when $x$ is a very general point. In \cite{Xu95} it is proved that the bound $m(m-1)<C^2$ holds when $x$ is a very general point.

\end{rem}

\subsection{Seshadri exceptional curves and elementary transformations.}\label{elemental}

An interesting property of Seshadri exceptional curves on ruled surfaces is that they are conserved by elementary transformations. We will denote by  $S'$ the elementary transformation of $S$ at $x\in S$ (see \cite{FuPe} for details). There is  a point $x'\in S'$ such that $S$ is the elementary transform of $S'$ at $x'$. Given a divisor $C$ on $S$, $C'$ will be denote its strict transform on $S'$. Let $e'$ be the invariant of $S'$ and $Y_0$ its minimum self intersection curve.

\begin{lemma}

Let $C\equiv aX_0+bf$ be a curve passing through $x$. If $a>0$ then 
$$
C'\equiv aY_0+(b-m+\frac{a}{2}(1-e+e'))f,\, \hbox{ and }\, mult_x'(C')=a-mult_x(C).
$$

\end{lemma}
{\bf Proof:} Suppose that $C$ passes through $x$ with multiplicity $m$. Let $C'\equiv a'Y_0+b'f$. Applying the properties of the elementary transformation (see \cite{FuPe}), it is clear that $a'=a$ and:
$$
C'^2=C^2+a^2-2am\impp 2a(b'-\frac{1}{2}ae')=2a(b-\frac{1}{2}ae)+a^2-2am.
$$
Moreover, $C'$ passes through $x'$ with multiplicity $m'=\tilde{f}\cdot \tilde{C}=f\cdot C-m=a-m$. \qed

\begin{prop}

Let $C\neq f$ be an irreducible curve passing through $x\in S$ and $C'$ its strict transform. Then $C$ is a Seshadri exceptional curve based at $x$ if and only if $C'\subset S'$ is a Seshadri exceptional curve based at $x'$. 

\end{prop}
{\bf Proof:} By the previous Lemma:
$$
C'^2-mult_{x'}(C')^2=C^2-mult_x(C)^2.
$$
This implies the equivalence between the Seshadri exceptional nature of $C$ and $C'$. \qed

From this result if we know the Seshadri exceptional curves at a point $x\in S$, we also know the Seshadri exceptional curves at $x'\in S'$. Let us see the relation between the influence areas of two of such curves $C$ and $C'$. A straightforward calculation shows:

\begin{lemma}

Let $C$ be a Seshadri exceptional curve based at $x$. Let $A\equiv aX_0+bf$ and $A'\equiv a'Y_0+b'f$. If
 $$\frac{a}{2(b-\frac{1}{2}ae)}=\lambda^2\hbox{ is a solution of }q_C(A,x)=\sqrt{A^2}$$
then
  $$\frac{a'}{2(b'-\frac{1}{2}a'e')}=\left({\frac{\lambda}{\lambda-1}}\right)^2\hbox{ is a solution of }q_{C'}(A',x')=\sqrt{A'^2}.
$$

\end{lemma} \qed

\begin{cor}\label{areas}

Let $C$ be a Seshadri exceptional curve based at $x$.  If
$$
Q_C=\{aX_0+bf\in Nef(S)|\, \lambda_1^2<\frac{a}{2(b-\frac{1}{2}ae)}<\lambda_2^2\}
$$
then 
$$
Q_{C'}=\left\{{a'Y_0+b'f\in Nef(S')|\, \left({\frac{\lambda_2}{\lambda_2-1}}\right)^2<\frac{a'}{2(b'-\frac{1}{2}a'e')}<\left({\frac{\lambda_1}{\lambda_1-1}}\right)^2}\right\}.
$$
\end{cor} \qed

Finally, by using elementary transformations we obtain a new bound on the multiplicity of a curve, which will be useful in some cases:

\begin{lemma}\label{cotasimple}

Let $S$ be a geometrically ruled surface and $C\equiv aX_0+bf$ with $a>1$ an irreducible curve passing through $x$ with multiplicity $m$. Then
$$
m\leq (b-\frac{1}{2}ae)+\frac{1}{2}a.
$$

\end{lemma}
{\bf Proof:} Let $S'$ be the elementary transformation of $S$ at $x$. We saw that:
$$
C'^2=C^2+a^2-2am.
$$
Moreover if $a>1$, $C'$ has non negative self-intersection. Thus:
$$
C^2+a^2-2am\geq 0 \impp 2a(b-\frac{1}{2}a)+a^2-2am\geq 0 \impp m\leq (b-\frac{1}{2}ae)+\frac{1}{2}a.
$$
\qed

\subsection{Seshadri constants on ruled surfaces with invariant $e>0$.}

\begin{lemma}
Let $S$ be a geometrically ruled surface  over a smooth curve $X$ with invariant $e>0$. The unique Seshadri-exceptional curves are the minimum self-intersection curve $X_0$ and the fibers $f$.
\end{lemma}
{\bf Proof:} Let $x\in X$ be an arbitrary point. Suppose that $C\equiv aX_0+bf$ is an irreducible curve passing through $x$ and different from $X_0$ and $f$. By the Lemma \ref{noestable}, we know that $a>0$ and $b>ae$. Applying the Corollary \ref{rank1} we deduce that $C$ is not a Seshadri exceptional curve. \qed

\begin{teo}\label{inestables}

Let $S$ be a geometrically ruled surface with invariant $e>0$. Let $A\equiv aX_0+bf$ be a nef linear system on $S$.  Then:

\begin{enumerate}

\item If $x\in X_0$ then $\epsilon(A,x)=min\{q_f(A,x),q_{X_0}(A,x)\}=min\{a,b-ae\}$.

\item If $x\not \in X_0$ then  $\epsilon(A,x)=q_f(A,x)=a$.

\end{enumerate}

\end{teo}
{\bf Proof:} It is a direct consequence of the previous Lemma. Note, that from the lemmas \ref{fibras} and  \ref{noestable}, we know that the influence area of $f$ fills the nef cone of $S$. \qed

\begin{rem}

For ruled surfaces with invariant $e>0$, we see that a nef linear system $A\equiv aX_0+bf$ reaches the expected value $\sqrt{A^2}$ only in two situations:

- when $e=1$ and $b=a$, at points $x\not \in X_0$.

- when $A\equiv b f$, at any point $x\in S$.

\end{rem}

\subsection{Seshadri constants on ruled surfaces with invariant $e\leq 0$.}

In this section we will work with a geometrically ruled surface $S$ over a smooth curve $X$ with invariant $e\leq 0$. It corresponds to a stable (or semi-stable) vector bundle of rank $2$. We have the following general bound for the Seshadri constants:

\begin{teo}\label{principal}

Let $S$ be a ruled surface with invariant $e\leq 0$. Let $A\equiv aX_0+bf$ be a nef divisor on $S$. Let $x$ be a point of $S$.

\begin{enumerate}

\item If  $e=0$ and $x$ lies on a curve numerically equivalent to $X_0$ then $$\epsilon(A,x)=min \{q_f(A,x), q_{X_0}(A,x)\}=min\{a,b\}.$$

\item In other case:

\begin{enumerate}

\item If $b-\frac{1}{2}ae\geq \frac{1}{2}a$ then
$$
\epsilon(A,x)=a.
$$

\item If $0\leq b-\frac{1}{2}ae\leq \frac{1}{2}a$ then
$$
2(b-\frac{1}{2}ae)\leq \epsilon(A,x)\leq \sqrt{A^2}=\sqrt{2a(b-\frac{1}{2}ae)}.
$$ 

\end{enumerate}

\end{enumerate}

\end{teo}
{\bf Proof:} Because $f$ is a Seshadri exceptional curve based at $x$, 
$$
\epsilon(A,x)\leq A\cdot f=a.
$$
If $e=0$ and $x$ lies in a curve $C\equiv X_0$, by the Corollary \ref{unisecantes} the unique Seshadri exceptional curves based at $x$ are $C$ and $f$. Moreover, 
$$
A^2=2ab=2q_{f}(A,x)q_{X_0}(A-x)\geq min\{A\cdot f, A\cdot X_0\}^2,
$$
so $\epsilon(A,x)=min\{A\cdot f, A\cdot X_0\}$.

If $e<0$, by the Lemma \ref{unicas}, we know that there are not unisecant Seshadri exceptional curves. Let $D\equiv kX_0+lf$ be an irreducible curve passing through $x$ with multiplicity $m$. We will apply the bound $m\leq \frac{1}{2}k+(l-\frac{1}{2}ke)$ obtained in the Lemma \ref{cotasimple}. We have:
$$
\dfrac{A\cdot D}{m}=\dfrac{a(l-\frac{1}{2}ke)+k(b-\frac{1}{2}ae)}{m}.
$$
If $b-\frac{1}{2}ae\geq \frac{1}{2}a$, we obtain:
$$
\dfrac{A\cdot D}{m}\geq \dfrac{a(l-\frac{1}{2}ke)+k(\frac{1}{2}a)}{m}\geq a
$$
and then 
$$
\epsilon(A,x)\geq a.
$$
If $b-\frac{1}{2}ae\leq \frac{1}{2}a$, we obtain:
$$
\dfrac{A\cdot D}{m}\geq \dfrac{2(b-\frac{1}{2}ae)(l-\frac{1}{2}ke)+k(b-\frac{1}{2}ae)}{m}\geq 2(b-\frac{1}{2}ae).
$$ \qed

\begin{cor}\label{trivial}

Let $S$ be the geometrically ruled surface $X\times \P^1$, where $X$ is a smooth curve. Let $A\equiv aX_0+bf$ be a nef divisor on $S$. Then for any point $x\in S$:
$$
\epsilon(A,x)=min\{q_{f}(A,x),q_{X_0}(A,x)\}=min\{a,b\}.
$$

\end{cor}
{\bf Proof:} When $S=X\times \P^1$ there is a curve $X\times \{t\}\sim X_0$ passing through any point. Thus, the result is a direct consequence of the previous Theorem. \qed

\begin{rem}

In this case the Seshadri constant does not depend on the point. Furthermore, a nef linear system $A\equiv aX_0+bf$ reaches the expected value $\sqrt{A^2}$ only when $a=2b$, $b=2a$, $a=0$ or $b=0$.

\end{rem}

We have seen that when $b-\frac{1}{2}ae\geq \frac{1}{2}a$ the Seshadri constant is determined by the Seshadri exceptional curve $f$. However, when $b-\frac{1}{2}ae<\frac{1}{2}a$ the fibres are not significant curves to compute the Seshadri constants. Let us see how in this range the problem is more delicate.

Note that when we have a linear system $A\equiv aX_0+bf$ with $\frac{1}{2}a$ near to $b-\frac{1}{2}ae$, the lower bound obtained in the Theorem \ref{principal} is near to the upper bound. In particular, if we take $a=2n$ and $b=ne+n-1$, then the bounds are:
$$
2(n-1)\leq \epsilon(A,x)\leq 2\sqrt{n(n-1)}
$$
with
$$
\lim_{n\rightarrow \infty} 2\sqrt{n(n-1)}-2(n-1)=1.
$$
In fact, we can construct ruled surfaces and linear systems where the Seshadri constant does not reach the upper bound, but it is as close as we wish.

\begin{teo}

Given a real number $\delta>0$ and a smooth curve $X$ of genus $>0$, there is a stable ruled surface $S$, an ample divisor $A$ on $S$ and a point $x\in S$ such that
$$
\sqrt{A^2}-\delta<\epsilon(A,x)<\sqrt{A^2}.
$$

\end{teo}
{\bf Proof:} Consider the ruled surface $S_0=\P(\Te_X\oplus \Te_X(\e))$ where $\e\in Pic^0(X)$ is a strict $n$-torsion point. Then there is a smooth curve $C\equiv nX_0$. This curve does not intersect $X_0$ and $X_1$. It is a Seshadri exceptional curve, because $C^2=0<1$. Let $y\in C$. If we make the elementary transformation of $S_o$ at $y$, we obtain a ruled surface $S$ with invariant $e=-1$. The strict transform of $C$ is a Seshadri exceptional curve $C'$ with a point $x$ of multiplicity $n-1$. Moreover, 
$$
C'^2=C^2+n^2-2n=n(n-2).
$$
Thus $C$ is an ample divisor on $S$ verifying:
$$
\epsilon(C',x)=\frac{C'^2}{m}=\frac{n(n-2)}{n-1}<\sqrt{C'^2}
$$
and:
$$
\lim_{n\rightarrow \infty} \sqrt{C'^2}-\frac{n(n-2)}{n-1}=\lim_{n\rightarrow \infty} \sqrt{n(n-2)}-\frac{n(n-2)}{n-1}=0.
$$ \qed

When $e=0$ we saw that the existence of curves $C\equiv X_0$ allows us to compute the Seshadri constants over their points. In general, the existence of curves $C$ with $C^2=0$ will give us valuable information to compute the constants. they are Seshadri exceptional curves at anyone of their points. By the Proposition \ref{estable} these curves can only appear on the linear systems:

-$|X_0|$, when $e=0$.

-$|nX_0+\frac{1}{2}ne|$ with $n\geq 2$.

Let us study the second case. Note that, by the Lemma \ref{cota}, an irreducible curve $C\equiv nX_0+\frac{1}{2}nef$ is smooth.

\begin{prop}\label{enesecantes}

Let $S$ be a ruled surface with invariant $e\leq 0$. Let $A\equiv aX_0+bf$ be a nef divisor and $C\equiv nX_0+\frac{1}{2}nef$ an irreducible curve with $n\geq 2$. Let $x\in C$:

\begin{enumerate}

\item If $b-\frac{1}{2}ae\geq \frac{1}{2}a$ then $$\epsilon(A,x)=a.$$

\item If $\frac{2}{n^2}a\leq b-\frac{1}{2}ae\leq \frac{1}{2}a$ then $$\frac{1}{n}a+2(1-\frac{1}{n})(b-\frac{1}{2}ae)\leq \epsilon(A,x)\leq \sqrt{A^2}=\sqrt{2a(b-\frac{1}{2}ae)}.$$

\item If $\frac{1}{n^2-2n+2}a\leq b-\frac{1}{2}ae\leq \frac{2}{n^2}a$ then $$\frac{1}{n}a+2(1-\frac{1}{n})(b-\frac{1}{2}ae)\frac{1}{n}\leq \epsilon(A,x)\leq n(b-\frac{1}{2}ae).$$

\item If $0\leq b-\frac{1}{2}ae\leq \frac{1}{n^2-2n+2}a$ then $$\epsilon(A,x)=n(b-\frac{1}{2}ae).$$

\end{enumerate} 

\end{prop}
{\bf Proof:} If $b-\frac{1}{2}ae\geq \frac{1}{2}a$ the result follows from the Theorem \ref{principal}. In other case, we have:
$$
\epsilon(A,x)\leq A\cdot C=n(b-\frac{1}{2}ae).
$$
This bound is smaller than the general bound $\sqrt{A^2}$ when $b-\frac{1}{2}ae\leq \frac{2}{n^2}a$.

Let $D\equiv kX_0+lf$ be a curve passing through $x$ with multiplicity $m>0$  and different from $C$. Then
$$m\leq D\cdot C=n(l-\frac{1}{2}ke)\iff \dfrac{l-\frac{1}{2}ke}{m}\geq \frac{1}{n}.$$

Applying the bound $k\geq 2m-2(l-\frac{1}{2}ke)$ obtained in the Lemma \ref{cotasimple}, we have:

$$
\begin{array}{rl}
{\dfrac{A\cdot D}{m}=}&{ \dfrac{a(l-\frac{1}{2}ke)+k(b-\frac{1}{2}ae)}{m}\geq  2(b-\frac{1}{2}ae)+(a-2(b-\frac{1}{2}ae))\dfrac{l-\frac{1}{2}ke}{m}\geq}\\
{}&{}\\
{}&{\geq 2(b-\frac{1}{2}ae)+(a-2(b-\frac{1}{2}ae))\frac{1}{n}=\frac{1}{n}a+2(1-\frac{1}{n})(b-\frac{1}{2}ae).}\\
\end{array}
$$
This bound is smaller than the given by $C$ when $a\leq (n^2-2n+2)(b-\frac{1}{2}ae)$. \qed

In particular, this result allow us to compute the Seshadri constants at points in smooth $2$-secant curves.

\begin{cor}\label{bisecantes}

Let $S$ be a ruled surface with invariant $e\leq 0$.
Let $A\equiv aX_0+bf$ be a nef divisor and $C\equiv 2X_0+ef$ an irreducible curve. Let $x\in C$. It holds:

\begin{enumerate}

\item If $b-\frac{1}{2}ae\geq \frac{1}{2}a$ then $\epsilon(A,x)=a$.

\item If $0\leq b-\frac{1}{2}ae\leq \frac{1}{2}a$ then $\epsilon(A,x)=2(b-\frac{1}{2}ae)$.

\end{enumerate} 

\end{cor}

\begin{example}\label{2torsion}

Let $X$ be a smooth curve of genus $g>0$ and $\e\in Pic^0(X)$ a $2$-torsion point. We consider the ruled surface  $S=\P(\Te_X\oplus \Te_X(\e))$.

The linear system $|2X_0|$ is a $1$-dimensional family of curves. The unique reducible curves are $2X_0$ and $2X_1$. Moreover, there is a smooth irreducible curve on $|2X_0|$ passing through any point $x\not\in X_0\cup X_1$. Applying the Theorem \ref{principal} and the Corollary \ref{bisecantes}, we can compute the Seshadri constant of any nef divisor at any point of $S$:

Let $A\equiv aX_0+bf$ be a nef divisor on $S$. It holds:

\begin{enumerate}

\item If $x\in X_0\cup X_1$ then $\epsilon(A,x)=min\{a,b\}$.

\item If $x\not\in X_0\cup X_1$ then:

\begin{enumerate}

\item If $b\geq \frac{1}{2}a$, $\epsilon(A,x)=a$.

\item If $0\leq q\leq \frac{1}{2}a$, $\epsilon(A,x)=2b$.

\end{enumerate}

\end{enumerate}

\end{example}

\section{Seshadri constants on rational ruled surfaces.}

Let $S$ be a  rational ruled surface. It is well known that $S$ is decomposable. In particular, $S=S_{e}=\P(\Te_{P^1}\oplus \Te_{P^1}(-e))$. Applying the Theorem \ref{inestables} and the Corollary \ref{trivial} we can compute the Seshadri constants over any  rational ruled surface:

\begin{teo}

Let $e\geq 0$ and let $S_e=\P(\Te_{P^1}\oplus \Te_{P^1}(-e))$ be a rational ruled surface. Let $A\equiv aX_0+bf$ be a nef divisor on $S_e$. It holds:

\begin{enumerate}

\item When $e=0$, $\epsilon(A,x)=min\{q_{f}(A,x),q_{X_0}(A,x)\}=min\{a,b\}$ for any $x\in S_e$.

\item When $e>0$:

\begin{enumerate}

\item If $x\in X_0$ then $\epsilon(A,x)=min\{q_{f}(A,x),q_{X_0}(A,x)\}=min\{a,b-ae\}$.

\item If $x\not\in X_0$ then $\epsilon(A,x)=q_f(A,x)=a$.

\end{enumerate}

\end{enumerate}

\end{teo}

\section{Seshadri constants on elliptic ruled surfaces.}\label{elipticas}

Let $X$ be a elliptic smooth curve. The classification of the geometrically ruled surfaces over $X$ is well known (see  \cite{At57},  \cite{FuPe00} or \cite{Ha}). There are the following possibilites:

\begin{enumerate}

\item Decomposable elliptic ruled surfaces $\P(\Te_X\oplus \Te_X(\e))$ with invariant $e\geq 0$. For each $e$ there is a one dimension family of such curves parameterized by $\e\in Pic^e(X)$. They will be denoted by $S_{\e}$.

\item The indecomposable elliptic ruled surface with invariant $e=0$ and $\e=0$. It will be denoted by $S_0$.

\item The indecomposable elliptic ruled surface with invariant $e=-1$. It will be denoted by $S_{-1}$.

\end{enumerate}

When $S_{\e}$ is a decomposable elliptic ruled surface with invariant $e> 0$ or $\e=0$, then the Theorem \ref{inestables} and the Corollary \ref{trivial} give the Seshadri constant of any nef divisor at any point. The interesting cases are $S_{\e}$ with $e=0$ and $\e\not\sim 0$, $S_0$ and $S_1$.

Let us remember some well known basic facts about these surfaces:

$S_{\e}$ has two unisecant disjoint curves $X_0$ and $X_1$ with self intersection zero. If $\e\in Pic^0(X)$ is a strict $k$-torsion point, then there is a smooth irreducible curve $B_k\sim kX_0$ passing through any point $x\not\in X_0\cup X_1$. On the contrary, if $\e\in Pic^0(X)$ is a non-torsion point, then there are not irreducible curves $B_k\sim kX_0$ for $k>1$.

$S_0$ has a unique unisecant curve $X_0$ with self intersection zero. Moreover, there are not irreducible curves $B_k\sim kX_0$ for $k>1$.

The surface $S_{-1}$ has a specially interesting geometry. This surface can be described as a quotient of $X\times X$ (see \cite{FuPe00} for details):
$$
S_{-1}=\frac{X\times X}{(P,Q)\sim (Q,P)}.
$$
Let $S'$ be the elementary transformation of $S_{-1}$ at $[P,Q]$. From the properties of the elementary transform and the description of $S_{-1}$ given at \cite{FuPe00}, we have

- If $P=Q$ then $S'$ is the indecomposable elliptic ruled surface with invariant $e=0$.

- If $P\neq Q$ then $S'$ is the decomposable ruled surface $\P(\Te_X\oplus \Te_X(P-Q))$.

\begin{defin}

We define the (possibly reducible) curves on $S$:
$$
\begin{array}{l}
	{T_1:=\{[P,P]\in S_{-1}\}.}\\
 {T_k:=\{[P,Q]\in S_{-1}|P-Q\hbox{ is a strict $k$-torsion point in $Pic^0(X)$}\},\hbox{ $k>1$.}}\\
 {T:=\cup_{k\geq 2} T_k.}\\
\end{array} 
$$

\end{defin}

\begin{rem}

In fact, it is esay to check that the reducible components of $T_k$, are the quotient of curves $\{(z,z+a)\in X\times X|\, z\in X\}$ where $a$ is a strict $k$-torsion point. When $k=2$, they correspond to the three well known $2$-secant irreducible curves on $S$ with self intersection $0$ (see \cite{FuPe00}). When $k>2$, they are numerically equivalent to $4X_0+2f$.

\end{rem}
With this notation we have the following theorem (see \cite{FuPe} and \cite{FuPe00}):

\begin{teo}\label{transformadaselipticas}

Let $S_{-1}$ be the elliptic ruled surface with invariant $e=1$. Let $S'$ be the elementary transform of $S$ at $x\in S$. Then:
\begin{enumerate}

\item If $x\in T_1$ then $S'$ is the indecomposable elliptic ruled surface $S_0$ with invariant $e=0$.

\item If $x\in T_k$, with $k>1$ then $S'$ is the decomposable ruled surface $S_{\e}=\P(\Te_X\oplus \Te_X(\e))$, where $\e\in Pic^0(X)$ is a strict $k$-torsion point.

\item If $x\not\in T\cup T_1$ then $S'$ is the decomposable ruled surface $S_{\e}=\P(\Te_X\oplus \Te_X(\e))$, where $\e\in Pic^0(X)$ is a non-torsion point.

\end{enumerate}
\end{teo}

We will use the following method to compute the Seshadri constants on these elliptic ruled surfaces. By using a suitable involution we will construct divisors on $S_{\e}$ and $S_0$ with high multiplicity at the fixed points by the involution. We will study the irreducible components of these divisors. They will be the Seshadri exceptional curves of $S_{\e}$ and $S_0$. Finally, using elementary transformations we will extend the results to the elliptic ruled surface $S_{-1}$.

\subsection{Construction of divisors with high multiplicity.}

If we fix a group structure on $X$ we can consider the symmetric involution:
$$
(-1):X\times X\lrw X\times X; \qquad (z_1,z_2)\lrw (-z_1,-z_2).
$$
It is clear that induces an involution on $S_{-1}$:
$$
(-1):S_{-1}\lrw S_{-1}; \qquad [z_1,z_2]\lrw [-z_1,-z_2].
$$
If $x\in S_{-1}$ is a fixed point by the involution, then there is an induced involution $(-1):S'\lrw S'$ on the elementary transform $S'$. Note that the curve $\{[z,-z]\in S_{-1},P\in X\}$ is a curve of fixed points. It corresponds to a fixed generator $f_0$ of $S_{-1}$. By Theorem \ref{transformadaselipticas} if we choose a suitable point of $f_0$, we have an involution for any elliptic ruled surface with invariant $e=0$\footnote{In \cite{LaTu82} it appears a study of the automorphisms of the elliptic ruled surfaces. However, the automorphism of $S_{\e}$ induced by the symmetric involution on the base curve $X$ has been forgotten. It is proven that an automorphism $\sigma\in Aut(X)$ lifts to $P(\E)$ if and only if $\sigma^*\E\cong \E\otimes L$ with $L\in Pic(X)$. When $\E=\Te_X\oplus \Te_X(\e)$ with $\e\in Pic^0(X)$ the possibility $\sigma=(-1)_X$ and $L=\Te_X(-\e)$ is not considered. }.

First, we consider the ruled surface $\pi:S_{\e}=\P(\Te_X\oplus \Te_X(\e))\lrw X$ with $\e\in Pic^0(X)$, $\e\not\sim 0$. Let us study the involution $(-1):S_{\e}\lrw S_{\e}$. It is clear that:
$$
(-1)(X_0)=X_1;\qquad
(-1)(X_1)=X_0;\qquad
\pi_*(-1)=(-1)_X.
$$
Thus there are four invariant generators on $S_{-\e}$. Let $P f$ be one of them. The involution $(-1)$ restricted to $P f$ has to fixed points not in $X_0\cup X_1$. Let $x_0$ be one of these fixed points. Let $|L|$ be a linear system invariant by the involution. Then $|L|$ can be decomposed in two eigenspaces of invariant divisors:
$$
\begin{array}{l}
{|L|^+=\{D\in |L|\,|(-1)D=D \hbox{ and } mult_{x_0}(D)\hbox { is even}\}.}\\
{|L|^-=\{D\in |L|\,|(-1)D=D \hbox{ and } mult_{x_0}(D)\hbox { is odd}\}.}\\
\end{array}
$$
The divisors of $|L|^+$ (respectively $|L|^-$) will be called {\it even} (respectively {\it odd}) divisors. 

Consider the linear system $|aX_0+aX_1+Pf|$. It is invariant by the involution. We know that (see \cite{FuPe}):
$$
h^0(\Te_{S_{\e}}(aX_0+aX_1+Pf))=2a+1.
$$
In particular, we can give the following system of divisors spanning $|L|$:
$$
\{nX_0+(2a-n)X_1+P_nf,\, n=0,1,\ldots,2a\}
$$
where $P_n$ is a point of $X$ verifying $aX_0+aX_1+Pf\sim nX_0+(2a-n)X_1+P_nf$.

Note that for any $n<a$, the pencils $$\langle nX_0+(2a-n)X_1+P_nf, (2a-n)X_0+nX_1+P_{2a-n}f\rangle$$ are invariant by the involution. However they are not fixed by the involution. Thus, they contain an even divisor and an odd divisor. Finally, the divisor $aX_0+aX_1+Pf$ is an odd divisor. We deduce that:
$$
h^0(\Te_{S_{\e}}(aX_0+aX_1+Pf))^+=a; \qquad h^0(\Te_{S_{\e}}(aX_0+aX_1+Pf))^-=a+1.
$$
Now, we can use this decomposition to get a divisor with high multiplicity at $x_0$:

\begin{prop}
Let $S_{\e}=\P(\Te_X\oplus \Te_X(\e))$ be a elliptic decomposable ruled surface, with $\e\in Pic^0(X)$, $\e\not\sim 0$. There is a divisor $C_n\equiv 2n(n+1)X_0+f$ on $S_{\e}$ such that $mult_{x_0}(C_n)\geq 2n+1$.

\end{prop}
{\bf Proof:} We consider the linear subsystem $|n(n+1)X_0+n(n+1)X_1+Pf|^-$. The number of conditions for a odd divisor to have a point of multiplicity $2n+1$ at $x$ is:
$$
2+4+\ldots+2n=n(n+1).
$$
But, we saw that $h^0(\Te_{S_{\e}}(n(n+1)X_0+n(n+1)X_1+Pf))^-=n(n+1)+1$. \qed

\begin{cor}\label{altad}
Let $S_{\e}=\P(\Te_X\oplus \Te_X(\e))$ be a elliptic decomposable ruled surface, with $\e\in Pic^0(X)$, $\e\not\sim 0$. Let $x$ be a point of $S_{\e}$ verifying $x\not\in X_0\cup X_1$. There is a divisor $C_n\equiv 2n(n+1)X_0+f$ on $S_{\e}$ such that $mult_x(C_n)\geq 2n+1$.
\end{cor}
{\bf Proof:} By propositions $3.2$ and $5.1$ in \cite{LaTu82}, we know that there is an automorphism $\sigma$ of $S_{\e}$ such that $\sigma(x_0)=x$. \qed

\begin{rem}

These divisors verify $C_n^2=4n(n+1)<(2n+1)^2=mult_x(C_n^2)$, so they are good candidates to contain a Seshadri exceptional curve. In fact, we will see that they are irreducible curves when $\e$ is generic. They allow us to compute the Seshadri constants at any point for each elliptic ruled surface.

\end{rem}

This construction can be easily extended to the indecomposable elliptic ruled surface $S_0$. We can study the involution $(-1):S_0\lrw S_0$. Now,
$$
(-1)(X_0)=X_0; \qquad \pi_*(-1)=(-1)_X.
$$
Let $P f$ be one of the four invariant generators. There are two fixed points in $P f$. One of them is $X_0\cap P f$ and the other one is $x_0\not\in X_0$. Moreover, we have invariant linear system $|2n(n+1)X_0+Pf|$ . We could study the dimension of the subspaces of even and odd divisors. However, we can give an alternative argument:

Note, that the ruled surfaces $S_{\e}$ are obtained from $S_{-1}$ applying elementary transformation at points $x\not\in T_1$. When $x\in T_1$, the elementary transform of $S_{-1}$ specializes to $S_0$. The two disjoint curves $X_0,X_1\subset S_{\e}$ degenerate to the curve $X_0\subset S_0$. In this way, the family of divisors $C_n\subset S_{\e}$ will provide a divisor $C_n\equiv 2n(n+1)X_0+f\subset S_0$ of multiplicity $\geq 2n+1$ at a point $x_0\not\in X_0$. Since, there is an automorphism $\sigma_x$ of $S_0$ moving $x_0$ to any other point $x\not\in X_0$ (see \cite{LaTu82}) we deduce:

\begin{prop}\label{altai}
Let $S_0$ be an elliptic indecomposable ruled surface with invariant $e=0$. Let $x$ be a point of $S_{\e}$ verifying $x\not\in X_0$. There is a divisor $C_n\equiv 2n(n+1)X_0+f$ on $S_{\e}$ such that $mult_x(C_n)\geq 2n+1$. \qed
\end{prop}

\begin{rem}\label{puntocomun}

Note, that the divisors $C_n\subset S_0$ with high multiplicity at $x_0$ are on the linear systems $|2n(n+1)X_0+Pf|$. Any linear system $|aX_0+(b P)f|$ with $b=0,1$ has a base point at $P f\cup X_0$ (see \cite{FuPe00}). Thus all of the irreducible components of $C_n$ passes through $P f\cap X_0$. This will be useful to prove their irreducibility.
\end{rem}

\subsection{Seshadri constants on $S_{\e}=\P(\Te_X\oplus \Te_X(\e))$ with $\e\in Pic^0(X)$ a non-torsion point.}

We will work on the elliptic ruled surface $S_{\e}=\P(\Te_X\oplus \Te_X(\e))$ with $\e\in Pic^0(X)$ a non-torsion point. In this case, we will see that the divisors $C_n$ constructed before are irreducible and with multiplicity $2n+1$ at points $x\not\in X_0\cup X_1$. From this, we will give all Seshadri exceptional curves at any point of $S_{\e}$. Note, that because $\e$ is a non-torsion potin, there are not irreducible curves $B\equiv kX_0$ for any $k>1$.

\begin{prop}\label{notorsion}

If $x\in S_{\e}$ is a point not in $X_0\cup X_1$ then
there is an irreducible curve $C_n\equiv 2n(n+1)X_0+f$ passing through $x\not\in X_0\cup X_1$ with multiplicity $2n+1$.

\end{prop}
{\bf Proof:} By the Proposition \ref{altad} we know that there is an effective divisor $C_n\equiv 2n(n+1)X_0+f$ passing through $x$ with multiplicity at least $2n+1$. We can suppose that $C_n$ is an odd divisor invariant by the involution $(-1)$. Thus it can be decomposed into irreducible components in the following way: $C_n=a_n(X_0+X_1)+b_nf+D_n$, where $D_n$ is an irreducible curve.  

If $C$ contains a fibre $f$, then $D_n\equiv (n(n+1)-a_n)(X_0+X_1)$, and $mult_x(D_n)\geq 2n$. But this is not possible because  there are not irreducible curves $B\equiv kX_0$ with $k>1$. We deduce that $b_n=0$, $D_n\equiv d_nX_0+f$ with $d_n\leq  2n(n+1)$ even and $mult_{x}(D_n)\geq 2n+1$. Moreover, applying the bound for the multiplicity obtained in the Lemma \ref{cota}, we see that $mult_{x}(D_n)=2n+1$. 

Now, we will prove that $D_n=C_n$. We will procced by induction. 

When $n=1$, $D_1\equiv d_1X_0+f$ with $d_1\leq 4$. The multiplicity of $D_1$ at $x$ is $3$. By Lemma \ref{cota} we know that:
$$
3(3-1)\leq 2(d_2-1) \impp d_1\geq 4.
$$
Thus $D_1=C_1$.

Suppose that $D_{n-1}=C_{n-1}$. Since $mult_{x}(D_n)>mult_{x}(C_{n-1})$, we know that $D_n\neq C_{n-1}$. Then:
$$
D_n\cdot C_{n-1}\geq mult_{x}(D_n)mult_{x}(C_{n-1})\impp d_n\geq 2n(n+1)-1.
$$
Because $d_n$ is even we conclude that $D_n=C_n$. \qed

With this notation, the fibre $f$ corresponds to the curve $C_0$. Thus we have a family $\{C_n\}$ of Seshadri exceptional curves based at $x\not\in X_0\cup X_1$. If we compute their influence area we obtain:

\begin{lemma}\label{influencian}
$$
Q_{C_n}=\{A\equiv aX_0+bf\in Nef(S)|\, n^2< \frac{a}{2b} < (n+1)^2\}.
$$ 
\end{lemma}

From this we deduce that they are the unique Seshadri exceptional curves based at $x\not\in X_0\cup X_1$. The Seshadri exceptional curves at points of $X_0$ and $X_1$ are well known from the general results (see Corollary \ref{unicas}).

\begin{teo}
Let $S=\P(\Te_X\oplus \Te_X(\e))$ be a decomposable ruled surface over a smooth elliptic curve $X$ with invariant $e=0$ and $\e\in Pic^0(X)$ a non-torsion point. Let $X$ be a point of $S$.

\begin{enumerate}

\item If $x\in X_i$, $i=1,2$, the unique Seshadri exceptional curves based at $x$ are 
$$
\begin{array}{l}
{\{f,\hbox { with } mult_x(f)=1\}.}\\
{\{X_i\hbox { with } mult_{x}(X_i)=1\}.}\\
\end{array}
$$

\item If $x\not\in X_0\cup X_1$, the unique Seshadri exceptional curves based at $x$ are 
$$
\begin{array}{l}
{\{f,\hbox { with } mult_x(f)=1\}.}\\
{\{C_n\equiv 2n(n+1)X_0+f,\hbox { with } mult_{x}(C_n)=2n+1\}_{n\geq 1}. }\\ 
\end{array}  
$$ \qed

\end{enumerate}

\end{teo}

Now, we are in a position to give the Seshadri constant of any ample divisor at any point of $S_{\e}$.

\begin{teo}

Let $S_{\e}=\P(\Te_X\oplus \Te_X(\e))$ be a decomposable ruled surface over a smooth elliptic curve $X$ with invariant $e=0$ and $\e\in Pic^0(X)$ a non-torsion point. Let $A\equiv aX_0+bf$ be an ample divisor on $S$. Let $x\in S$. Then:

\begin{enumerate}

\item If $x\in X_0\cap X_1$ then $\epsilon(A,x)=min\{q_f(A,x),q_{X_0}(A,x)\}=min\{a,b\}$.

\item If $x\not\in X_0\cap X_1$ then
$$
\epsilon(A,x)=q_{C_n}(A,x)=\frac{2n(n+1)b+a}{2n+1},
\hbox{ with }
n=\left[\sqrt{\frac{a}{2b}}\right].
$$
\end{enumerate}

\end{teo}

\subsection{Seshadri constants on $S=\P(\Te_X\oplus \Te_X(\e))$ with $\e\in Pic^0(X)$ a strict $k$-torsion point.}

We will work on the elliptic ruled surface $S_{\e}=\P(\Te_X\oplus \Te_X(\e))$ with $\e\in Pic^0(X)$ a strict $k$-torsion point. The main different with the previous case is the existence of irreducible smooth curves $B_k\sim kX_0$ passing through any point $x\not\in X_0\cup X_1$.

\begin{prop}

If $x\in S_0$ is a point not in $X_0$ and
If $n<\frac{k}{2}$ then there is an irreducible curve $C_n\equiv 2n(n+1)X_0+f$ passing through $x$ with multiplicity $2n+1$.

\end{prop}
{\bf Proof:} We can apply the arguments of the Proposition \ref{notorsion} if we prove that $B_k$ is not a reducible component of $C_n$. When $n=1$ it is clear that $B_k$ and $f$  are not reducible components of $C_1\equiv 4X_0+f$. It is sufficient to use that $mult_{x}(C_1)\geq 3$. Thus $C_1$ is irreducible and by the Lemma \ref{cota} the multiplicity is exactly $3$.

Suppose that the theorem holds for $C_{n-1}$. If $C_{n-1}$ is a reducible component of $C_n$, then $mult_{x}(C_n-C_{n-1})\geq 2n+1-(2n-1)=2$. From this $2B_k$ must be a reducible component of $C_{n}-C_{n-1}$. Then $2n(n+1)-2(n-1)n\geq 2k$, so $k\leq 2n$. This contradicts the hypothesis. 

If $B_k$ is a reducible component of $C_n$ then $C_n-B_k$ is an effective divisor passing through $x$ with multiplicity $\geq 2n$. Thus:
$$
C_{n-1}\cdot (C_n-B_k)\geq (2n-1)2n\impp 2n(n+1)\geq k.
$$
But this contradicts the hypothesis again. \qed

Now, we compute the influence area of $B_k$.

\begin{lemma}\label{areak}

$$Q_{B_k}=\{A\sim aX_0+bf\in Nef(S)|\,
\frac{a}{2b} > \frac{k^2}{4}\}. 
$$ \qed

\end{lemma}

The influence area of the curves $C_n$ was computed at the Lemma \ref{influencian}. Thus, we know the set of Seshadri exceptional curves based at any point of $S_{\e}$ and we can compute the Seshadri constants.

Note that when $k$ is odd and  $n=\frac{k-1}{2}$, the influence are of $B_k$ and $C_{n}$  intersect. In this case we can check that if $\frac{a}{2b}\leq \frac{k^2+1}{4}$ then $q_{C_n}(A,x)\leq q_{B_k}(A,x)$ for any $A\equiv aX_0+bf$.

\begin{teo}

Let $S=\P(\Te_X\oplus \Te_X(\e))$ be a decomposable ruled surface over a smooth elliptic curve $X$ with invariant $e=0$ and $\e\in Pic^0(X)$ a $k$-torsion point. Let $x$ be a point of $S$.
\begin{enumerate}

\item If $x\in X_i$, $i=1,2$, the unique Seshadri exceptional curves based at $x$ are 
$$
\begin{array}{l}
{\{f,\hbox { with } mult_x(f)=1\}.}\\
{\{X_i\hbox { with } mult_{x}(X_i)=1\}.}\\
\end{array}
$$

\item If $x\not\in X_0\cup X_1$, the unique Seshadri exceptional curves based at $x$ are 
$$
\begin{array}{l}
{\{f,\hbox { with } mult_x(f)=1\}.}\\
{\{B_k\equiv kX_0,\hbox { with } mult_x(B_k)=1\}.}\\
{\{C_n\equiv 2n(n+1)X_0+f,\hbox { with } mult_{x}(C_n)=2n+1\}_{1\leq n\leq \left[\frac{k-1}{2}\right]}.}\\
\end{array}
$$

\end{enumerate}

\end{teo}

\begin{teo}

Let $S=\P(\Te_X\oplus \Te_X(\e))$ be a decomposable ruled surface over a smooth elliptic curve $X$ with invariant $e=0$ and $\e\in Pic^0(X)$ a $k$-torsion point. Let $A\equiv aX_0+bf$ be an ample divisor on $S$. Let $x\in S$. Then:

\begin{enumerate}

\item If $x\in X_0\cap X_1$ then $\epsilon(A,x)=min\{q_f(A,x),q_{X_0}(A,x)\}=min\{a,b\}$.

\item If $x\not\in X_0\cap X_1$ and

\begin{enumerate}

\item $\frac{a}{2b}\leq \frac{k^2}{4}$ when $k$ is even or $\frac{a}{2b}\leq \frac{k^2+1}{4}$ when $k$ is odd, then

$$
\epsilon(A,x)=q_{C_n}(A,x)=\frac{2n(n+1)b+a}{2n+1},
\hbox{ with }
n=\left[\sqrt{\frac{a}{2b}}\right].
$$

\item $\frac{a}{2b}\geq \frac{k^2}{4}$ when $k$ is even or $\frac{a}{2b}\geq \frac{k^2+1}{4}$ when $k$ is odd, then

$$
\epsilon(A,x)=q_{B_k}(A,x)=kb.
$$

\end{enumerate}

\end{enumerate}

\end{teo}

\subsection{Seshadri constants on the indecomposable elliptic ru\-led surface $S_0$ with invariant $e=0$.}

We will work on the elliptic indecomposable ruled surface $S_0$ with invariant $e=0$. This case is similar to the decomposable ruled surface $S_{\e}$ where $\e\in Pic^0(X)$ is a non-torsion point. The divisors $C_n\equiv 2n(n+1)X_0+f$ passing through $x\not\in X_0$ are irreducible and have multiplicity $2n+1$ at $x$. 

To check this, we can apply most of the arguments of the proof of the Proposition \ref{notorsion}. The difference is that $d_n$ could be odd. However, we know that the curves $D_n$ has a common point at $X_0\cap Pf$ for any $n\geq 1$ (see Remark \ref{puntocomun}). Thus the condition:
$$
D_n\cdot C_{n-1}\geq mult_{x}(D_n)mult_{x}(C_{n-1})
$$
in the proof of the Proposition \ref{notorsion}, can be changed to:
$$
D_n\cdot C_{n-1}\geq mult_{x}(D_n)mult_{x}(C_{n-1})+1.
$$
From this, we obtain $d_n\geq 2n(n+1)$ and then $C_n=D_n$. \qed
\begin{prop}
If $x\in S_0$ is a point not in $X_0$ then
there is an irreducible curve $C_n\equiv 2n(n+1)X_0+f$ passing through $x\not\in X_0$ with multiplicity $2n+1$.

\end{prop}

Now, if we replace the curve $X_1$ by $X_0$, the results that we saw for $S_{\e}$ with $\e$ generic are valid in this case. We obtain:

\begin{teo}

Let $S_0$ be the indecomposable ruled surface over a smooth elliptic curve $X$ with invariant $e=0$. Let $A\equiv aX_0+bf$ be an ample divisor on $S$. Let $x\in S_0$. Then:

\begin{enumerate}

\item If $x\in X_0$ then $\epsilon(A,x)=min\{q_f(A,x),q_{X_0}(A,x)\}=min\{a,b\}$.

\item If $x\not\in X_0$ then
$$
\epsilon(A,x)=q_{C_n}(A,x)=\frac{2n(n+1)b+a}{2n+1},
\hbox{ with }
n=\left[\sqrt{\frac{a}{2b}}\right].
$$
\end{enumerate}

\end{teo}

\subsection{Seshadri constants on the elliptic ruled surface with invariant $e=-1$.}

Let $S_{-1}$ be the indecomposable ruled surface with invariant $e=-1$. Let $x$ be a point of $S_{-1}$.  If we make the elementary transformation of $S_{-1}$ at $x\in S_{-1}$, we obtain a ruled surface $S'$ with invariant $e=0$. The Theorem \ref{transformadaselipticas} describes with precision this ruled surface $S'$ depending on the position of the point $x$. 

We have computed al the Seshadri exceptional curves of the elliptic ruled surfaces with invariant $e=0$. Applying the results of the Section \ref{elemental} we know that their strict transforms will be the Seshadri exceptional curves based at $x$. In particular, we have:
$$
\begin{array}{lcl}
{C'_n\equiv 2n(n+1)X_0-2nf}&{ \qquad \hbox{with}\qquad}&{ mult_x(C'_n)=2n^2-1.}\\
{B'_k\equiv kX_0-f}&{ \qquad \hbox{with}\qquad}&{ mult_x(B'_k)=k-1.}\\
\end{array}
$$

\begin{teo}
Let $S_{-1}$ be an elliptic ruled surface with invariant $e=-1$. Let $x$ be a point of $S_{-1}$:

\begin{enumerate}

\item If $x\in T_k$ with $k>1$, then the unique Seshadri exceptional curves based at $x$ are 
$$
\begin{array}{l}
{\{f,\hbox { with } mult_x(f)=1\}.}\\
{\{B'_k\equiv kX_0-f,\hbox { with } mult_x(B'_k)=0\}.}\\
{\{C'_n\equiv 2n(n+1)X_0-2nf,\hbox { wih } mult_{x}(C'_n)=2n^2-1\}_{1\leq n\leq \left[\frac{k-1}{2}\right]}.}\\
\end{array}
$$ 

\item  If $x\not\in T$ then the unique Seshadri exceptional curves based at $x$ are 
$$
\begin{array}{l}
{\{f,\hbox { with } mult_x(f)=1\}.}\\
{\{C'_n\equiv 2n(n+1)X_0-2nf,\hbox { with } mult_{x}(C'_n)=2n^2-1\}_{n\geq 1}.}\\
\end{array}
$$ \qed

\end{enumerate}

\end{teo}

From the Corollary \ref{areas} and the Lemmas \ref{influencian} and \ref{areak}, we obtain the influence areas of these exceptional Seshadri curves:

\begin{lemma}

Let $A\equiv aX_0+bf$ be an ample divisor. Then:
$$
\begin{array}{rl}
{\hbox{(i)}}&{Q_{C'n}=\{A\equiv aX_0+bf\in Nef(S)|\, 
\left({1+\frac{1}{n}}\right)^2< \frac{a}{2(b+\frac{1}{2}a)} < \left({1+\frac{1}{n-1}}\right)^2\}}\\
{\hbox{(ii)}}&{Q_{B'_k}=\{A\equiv aX_0+bf\in Nef(S)|\,
\left({1+\frac{1}{\frac{k}{2}-1}}\right)^2< \frac{a}{2(b+\frac{1}{2}a)}\}}\\
\end{array}
$$

\end{lemma} \qed

Finally, with this information, we can compute the Seshadri constant of any ample divisor on $S_{-1}$:

\begin{teo}

Let $S_{-1}$ be an elliptic ruled surface with invariant $e=-1$. Let $A\equiv aX_0+bf$ be an ample divisor on $S_{-1}$. Let $x\in S_{-1}$.

\begin{enumerate}

\item If $x\in T_k$ with $k>1$ and 

\begin{enumerate}

\item $\frac{a}{2(b+\frac{1}{2}a)}\leq 1$ then
$$
\epsilon(A,x)=q_f(A,x)=a.
$$

\item $1<\frac{a}{2(b+\frac{1}{2}a)}\leq \frac{k^2}{(k-2)^2}$ when $k$ is even or $1<\frac{a}{2(b+\frac{1}{2}a)}\leq \frac{k^2+1}{k^2-4k+5}$ when $k$ is odd, then

$$
\epsilon(A,x)=q_{B'_k}(A,x)=\frac{(k-1)a+kb}{k-1}.
$$

\item $\frac{a}{2(b+\frac{1}{2}a)}\geq \frac{k^2}{(k-2)^2}$ when $k$ is even or $\frac{a}{2(b+\frac{1}{2}a)}\geq \frac{k^2+1}{k^2-4k+5}$ when $k$ is odd, then

$$
\epsilon(A,x)=q_{C'_n}(A,x)=\frac{2n(na+(n+1)b)}{2n^2-1}
,\;\;\;
n=\left[\frac{1}{\sqrt{\frac{a}{2(b+\frac{1}{2}a)}}-1}\right]+1.
$$

\end{enumerate}

\item If $x\not\in T$ and 
\begin{enumerate}

\item $\frac{a}{2(b+\frac{1}{2}a)}\leq 1$ then
$$
\epsilon(A,x)=q_{f}(A,x)=a.
$$

\item $\frac{a}{2(b+\frac{1}{2}a)}> 1$ then
$$
\epsilon(A,x)=q_{C'_n}(A,x)=\frac{2n(na+(n+1)b)}{2n^2-1}
,\;\;\;
n=\left[\frac{1}{\sqrt{\frac{a}{2(b+\frac{1}{2}a)}}-1}\right]+1.
$$

\end{enumerate}

\end{enumerate}

\end{teo}

{\bf E-mail:} lfuentes@udc.es

Luis Fuentes Garc\'{\i}a. 

Departamento de M\'etodos Matem\'aticos y Representaci\'on.

E.T.S. de Ingenieros de Caminos, Canales y Puertos. 

Universidad de A Coru\~na. Campus de Elvi\~na. 15192 A Coru\~na (SPAIN)

\end{document}